\documentclass[12pt,a4paper]{amsart}
\usepackage{amsmath, amssymb, bm}
\usepackage[all]{xy}

\def \LL {\skp \setlength{\leftskip}{8pt}}  
\def \LB {\setlength{\leftskip}{0pt}}   
\def \sst {\scriptstyle}

\def \c {\colon}

\def \q {\qquad}
\def \qq {\qquad \qquad}
\def \qqq {\qquad \qquad \qquad \qquad}

\def \sep {\,  |  \,}     
\def \bu {{\scriptscriptstyle\bullet}}
\def \skp {\medskip}
\def \ndt {\noindent}
\def \Ndt {\medskip  \noindent}
\def \tilde {\raise.17ex\hbox{$\scriptstyle\mathtt{\sim}$}}   

\def \adj {\dashv}
\def \dt {\raisebox{-.35ex}{${\mbox{\fontsize{18}{10}\selectfont\ensuremath{\cdot}}}$}}  
\def \todot {{\; \dt \hspace{-2.2ex} \to \,}}   

\def \ard {\ar@{-->}}
\def \arp {\ar@{.>}}
\def \are {\ar@{->>}}
\def \aru {\ar@{=}}
\def \arv {\ar@{}}   
\def \arl {\ar@{-}}    
\def \arld {\ar@{--}}    
\def \arlp {\ar@{..}}    

\def \ti {\! \times \!}
\def \jo {{\, {\scriptstyle{\vee}} \,}}
\def \me {{\, {\scriptstyle{\wedge}} \, }}
\def \we {{{\scriptstyle{\wedge}}}}    
\def \ci {{\raise.3ex\hbox{$ \scriptscriptstyle\circ $}}}
\def \Pro {\raisebox{0.45ex}{${\mbox{\fontsize{10}{10}\selectfont\ensuremath{\prod}}}$}}

\def \Cup {\raisebox{0.45ex}{${\mbox{\fontsize{9}{10}\selectfont\ensuremath{\bigcup}}}$}}
\def \Cap {\raisebox{0.45ex}{${\mbox{\fontsize{9}{10}\selectfont\ensuremath{\bigcap}}}$}}
\def \setm {{\raise.4ex\hbox{$ \, \scriptscriptstyle{\setminus} \; $}}}
\def \sub {\subset}
\def \sups {\supset}
\def\le{\leqslant}
\def\ge{\geqslant}

\def \eq {\! \sim \!}

\def \and {\mbox{ and }}

\def \int {{\rm int}}
\def \Ker {{\rm Ker}}
\def \Coker {{\rm Coker}}
\def \Im {{\rm Im}}

\def\dd {\partial}
\def\ddm {\partial^-}
\def\ddp {\partial^+}

\def \al {\alpha}

\def \ga {\gamma}
\def \ep {\varepsilon}
\def \ka {\kappa}
\def \si {\sigma}
\let \ph \varphi
\let \om \omega
\def \De {\Delta}
\def \Si {\Sigma}
\def \Om {\Omega}

\def \Set {\mathsf{Set}}
\def \Top {\mathsf{Top}}

\def \es {\varnothing}  
\def \sing  {\{*\}}
\def \bbN {{\mathbb{N}}}  
\def \bbZ {\mathbb{Z}}
\def \bbQ {\mathbb{Q}}
\def \bbI {{\mathbb{I}}}   
\def \bbR {{\mathbb{R}}}
\def \bbS {\mathbb{S}}
\def \bbD {\mathbb{D}}

\begin{document}

\title[Smash product and regular associativity]{A note on the smash product \\ and regular associativity}

\author[M. Grandis]{Marco Grandis}

 \address{Marco Grandis, Dipartimento di Matematica, Universit\`a di Genova, 16146-Genova, Italy.}
 \email{grandis@dima.unige.it}

 \subjclass{55Pxx, 18D15}

\keywords{Smash product, homotopy constructions, colax monoidal category}

\begin{abstract}

We want to study the smash product of pointed topological spaces, in an organic way and full 
generality, without relying on some `convenient subcategory'. The $n$-ary smash product has a 
`colax' form of associativity, which supplies a categorical framework for the properties of this 
operation and its connection with the function spaces.

	Various concrete computations of smash products are given, including a large class of cases 
where associativity fails.

	Lax and colax monoidal structures are unusual and interesting, in category theory. Some parts of 
this note will be obvious to a topologist and others to a categorist, in order to take into account both 
backgrounds.
\end{abstract}

 \maketitle 

\section*{Introduction}\label{Intro}

	The smash product $X \me Y$ of pointed topological spaces is not associative, unless in 
a weak form discussed below. Moreover a pointed space is exponentiable (for the smash product) 
if and only if the underlying topological space is exponentiable (for the cartesian product), that is core 
compact -- a kind of local compactness. These facts are closely related: exponentiable pointed 
spaces give associative triples for the smash product (but they are not closed under this operation).

	Because of these reasons the smash product is generally studied in a suitable subcategory of 
(pointed) spaces, like CW-complexes or (some version of) compactly generated spaces. Or is left out, 
as in many textbooks on Algebraic Topology.

\skp	Here we want to show that the category $\Top_\bu$ of pointed spaces need not be viewed as 
a hopeless structure that should be replaced by a convenient one: it is a {\em colax} monoidal category, 
a structure recently introduced in category theory \cite {Ln}, where partial associativity properties can 
be organically investigated. This moderate interplay of topology and category theory might hopefully 
be of interest in both fields, although this will require working out issues that can be obvious in one 
of these domains, or the other.

	More precisely, the smash product can be extended to an $n$-ary operation $ X_1 \me ... \me X_n$ 
with a natural system of associativity comparisons, generally non-invertible; this system is easily proved 
to be coherent, making $\Top_\bu$ into a colax monoidal category (see \ref{1.7}).

	In particular, a triple of pointed spaces has two comparison maps
    \begin{equation} \begin{array}{c}
\ga'\c X \me Y \me Z \to X \me (Y \me Z),   \q   \ga'(x \me y \me z)  =  x \me (y \me z),
\\[5pt]
\ga''\c X \me Y \me Z \to (X \me Y) \me Z,   \q   \ga''(x \me y \me z)  =  (x \me y) \me z,
    \label{01} \end{array} \end{equation}
and we say that a ternary smash product $P' = X \me (Y \me Z)$ or $P'' = (X \me Y) \me Z$ is 
{\em regular} if the corresponding comparison is invertible (see \ref{1.8}), which means that $P'$ or 
$P''$ is a quotient of the pointed space $X \me Y \me Z$ (namely, the quotient that collapses the 
subspace $X \jo Y \jo Z$ to the basepoint). 

	If both are, $X \me (Y \me Z)$ and $(X \me Y) \me Z$ coincide with $X \me Y \me Z$, up to 
canonical homeomorphism, and we say that the triple $(X, Y, Z)$ is {\em regularly associative}. This 
makes sense in any lax or colax monoidal category, and is stronger than mere associativity, 
as defined in \ref{1.3}: see \ref{1.8}(c). This property is already considered in an article by Carboni and 
Janelidze \cite{CrJ}, in a general context of pointed categories, with a different goal: studying conditions 
under which the property universally holds (see the Note \ref{1.8}(e)). It should be interesting to develop 
a common extension of their results and the present ones.

\Ndt {\em Outline}. The first two sections analyse the partial associativity properties of the smash 
product and their links with the exponentiability of pointed spaces. Our main contributions are 
concerned with the colax structure of the smash product (in \ref{1.7}) and the property of regular 
associativity (in \ref{1.8}); the latter is granted by conditions of exponentiability or compactness of a pair 
of the spaces involved (Theorem \ref{2.2}).

	Section 3 reviews in this light some facts, mostly known in topology, on the homotopy constructions 
related to the cofibre and fibre sequences of a map: pointed cylinders, cones and suspensions can be 
expressed as smash products, and commute with each other; their adjoints -- path spaces, cocones and 
loop spaces -- are the corresponding exponentials.

	Finally, Section 4 explores a large class of triples $(\bbQ, Y, X)$ of pointed spaces where 
associativity fails (Theorem \ref{4.8}), including all cases where $X$ and $Y$ are non-degenerate 
intervals and $X$ is not compact; it also contains the well-known triple $(\bbQ, \bbQ, \bbN)$ proposed 
by Puppe; the proof is shorter and simpler than that reported in \cite{MyS}, Section 1.7, for this case. 
The smash product of euclidean intervals is studied in the first part of the section, after the 
`easy' compact case, in \ref{1.5}.

\Ndt {\em Literature}. The smash product and exponential law of {\em general} pointed spaces is rarely 
considered. One can see Maunder \cite{Mn}, Section 6.2, tom Dieck \cite{Di} and also Hatcher \cite{Ha} 
(where associativity is left out). The general theory of lax monoidal categories can be found in 
Leinster's book \cite{Ln}.

\vskip 5pt

	In a different approach to the exponential law in topology, several `convenient categories'  have 
been proposed and studied, aiming to replace the category of topological spaces by some full 
(possibly coreflective) subcategory which is cartesian closed -- at the price of modifying the cartesian 
product and subspaces (at least). Pointed spaces are similarly replaced.

	An article by Escard\'o, Lawson and Simpson \cite{EsLS} studies and compares various candidates, 
like compactly generated Hausdorff spaces (also known as Kelley spaces), compactly generated spaces, 
locally compactly generated spaces, core compactly generated spaces and sequential spaces -- in the 
terminology of this article. Other references on this domain can be found there.

	The smash product in (different) `convenient categories' of pointed topological spaces is studied in 
various books and articles, like \cite{Wh, My, MyS, AgGP, ShYH}; the last is based on $\De$-generated spaces, 
called D-spaces in the domain of diffeology \cite{IZ}.

\vskip 5pt

	Finally, the article \cite{CrJ} mentioned above deals with the smash product of pointed objects, 
constructed on a ground category. Under suitable conditions on the latter this operation is proved 
to be associative and distributive with respect to finite sums. There are several examples, including -- 
obviously -- the category of compactly generated spaces and the opposite of the category of 
commutative algebras.

\Ndt {\em Notation and terminology}. $\Top$ and $\Top_\bu$ denote the categories of (topological) 
spaces and pointed spaces. A {\em map} is a continuous mapping, and {\em nbd} means `neighbourhood'. 
A compact space is not assumed to be Hausdorff. A space is locally compact if every point has a 
fundamental system of compact nbds.

	The symbols $ \bbR^n $ and $ \bbS^n $ denote the usual euclidean spaces. $ \bbI $ is the standard 
euclidean interval $ [0, 1]$; open and semi-open real intervals are also written with square brackets, 
like $ ]0, 1[ $ and $ [0, 1[$. The symbol $ \sub $ denotes weak inclusion. The binary variable $\al$ takes 
values  0, 1, also written as $-, +$ in superscripts and subscripts. Marginal remarks are written in small 
characters.

\Ndt {\em Acknowledgements}. I wish to express my gratitude to Ettore Carletti, George Janelidze and 
Francesca Cagliari for many discussions, helpful suggestions and invaluable help in providing references.

\section{Smash product of pointed spaces}\label{s1}

		After two preliminary subsections we work in the category $ \Top_\bu $ of pointed spaces 
and pointed maps. The basepoint of a pointed space $ X $ is usually written as $ 0_X$, or 0, and 
the underlying topological space is written as $ |X|$. The notation as a pair $ (|X|, 0_X) $ is rarely used. 
A map between pointed spaces is understood to be pointed, by default.

	The smash product $ X \me Y $ of pointed spaces -- a sort of tensor product related to `bipointed' 
maps -- is associative under suitable assumptions, studied here and in the next section. Formally it is a 
`colax' monoidal product, see \ref{1.7}. Concretely, we also begin to study the smash product of the 
euclidean intervals, in \ref{1.5}.

\subsection{Exponentiable spaces and function spaces}\label{1.1}
(a) We recall that a topological space $ A $ is {\em exponentiable} (for the cartesian product) if the functor 
$ - \ti A\c \Top \to \Top $ has a right adjoint, written as $ (-)^A\c \Top \to \Top$.

	We may assume that the {\em exponential space} $\; Y^A \;$ is the hom-set $ \Top(A, Y) $ 
equipped with the {\em exponential topology}, strictly determined by the adjunction. The adjunction is 
then expressed by a natural bijection in the variables $ X, Y$
    \begin{equation} \begin{array}{c}
\ph_{XY}\c \Top(X \ti A, Y) \to \Top(X, Y^A),
\\[5pt]
(f\c X \ti A \to Y)  \mapsto  (g\c X \to Y^A),   \;\;\;   g(x)  =  f(x, -)\c A \to Y,
    \label{1.1.1} \end{array} \end{equation}
corresponding to the exponential law of $\Set$.

\Ndt (b) The exponential topology $ Y^A $ is characterised as the finest topology of $ \Top(A, Y) $ 
that makes all transpose mappings $ g\c X \to Y^A $ continuous, for every topological space $ X $ 
and every map $ f\c X \ti A \to Y $ (\cite{EsLS}, Proposition 5.13).

\Ndt (c) If $ A, B $ are exponentiable spaces, also $ A \ti B $ is exponentiable and there is a canonical 
isomorphism of functors
    \begin{equation}
( - )^{A \ti B}  =  (( - )^B)^A,
    \label{1.1.2} \end{equation}
an obvious consequence of the composition of right adjoints (to the functors $ - \ti A $ and $ - \ti B$). We 
also note that the endofunctors $ ( - )^A $ and $ ( - )^B $ commute, as their left adjoints.

\subsection{Complements and literature}\label{1.2}
(a) Every locally compact space $ A $ (in the sense that each point has a local basis of compact nbds) 
is exponentiable. In this case, the exponential $Y^A$ is the function set $\Top(A, Y)$ with the 
compact-open topology. For instance, this is proved in \cite{Di}, Propositions 2.4.1 and 2.4.3.

\Ndt (b) As a partial converse, every exponentiable space which is Hausdorff must be locally compact, 
as proved in \cite{He}, 7.3.18. In particular, the rational line $ \bbQ $ is not exponentiable; this is also 
proved here, in \ref{4.8}. Let us recall that a compact non-Hausdorff space need not be locally compact; 
for instance, the one-point compactification of $\bbQ$ is not.

\Ndt (c) Exponentiable spaces are characterised as the {\em core compact} ones, which means that 
every open nbd $ V $ of a point $ x $ contains an open nbd $U$ such that every open cover of $ V $ 
has a finite subcover of $ U$.

	This characterisation is essentially due to Day and Kelly \cite{DaK}; for a recent exposition 
see \cite{EsH}.

\Ndt (d) As another characterisation, the space $ A $ is exponentiable if and only if the functor 
$ - \ti A $ preserves quotient maps (that is surjective maps whose codomain has the final topology for 
the map itself).

\LL \begin{small}

\vskip -2pt

	The necessity of this condition is obvious, because topological quotients are particular colimits. 
As to sufficiency, any functor $ - \ti A $ preserves sums, and the rest can be deduced from the 
Adjoint Functor Theorem \cite{Ma}: the `solution set condition' for the existence of the right adjoint is 
easily verified.

\end{small} \LB

\Ndt (e) A {\em bug-eyed interval}. On the other hand, a compact, locally compact space (is 
exponentiable and) need not be Hausdorff. One can easily form an example as a compact subspace 
$ A $ of the {\em bug-eyed line}, the usual instance of a locally euclidean space which is not Hausdorff.

\LL \begin{small}

\vskip -2pt

	The latter is the quotient $ L = (\bbR \ti \{0, 1\})/R$, modulo the equivalence relation that identifies 
all pairs $ (x, 0) $ and $ (x, 1) $ for $ x \neq 0$, leaving two points $ 0_\al = [(0, \al)] $ (the {\em bug eyes}) 
with no disjoint pair of nbds (for $ \al = 0, 1$). The space $ L $ is locally compact, being locally euclidean. 
The image $ A $ of $ \bbI \times \{0, 1\} $ in $ L $ is a compact, locally compact, non-Hausdorff subspace.

\end{small} \LB

\subsection{Smash product of pointed spaces}\label{1.3}
(a) We recall that $ \Top_\bu $ is a pointed category, with zero object the singleton $ \sing$, and zero maps 
$ X \to \sing \to Y$. Every map $ f\c X \to Y $ has a kernel $ \Ker f $ (the preimage of the basepoint of 
$ Y$) and a cokernel $ \Coker f $ (the quotient $ Y/f(X) $ that collapses the image of $ f $ to the new 
basepoint).

\Ndt (b) We are interested in an operation, the {\em smash product} $ X \me Y$, or {\em reduced product}
    \begin{equation}
- \me -\c \Top_\bu \ti \Top_\bu \to \Top_\bu,   \q\;   X \me Y  =  (X \ti Y) / (X \jo Y),
    \label{1.3.1} \end{equation}
which is symmetric but not associative (unless in a colax form to be investigated later).

\vskip 2pt

	Here the {\em wedge} $X \jo Y  =  X \ti \{0\} \cup \{0\} \ti Y $ is the pointed sum (that is the categorical 
sum of pointed spaces), realised as a subspace of $X \ti Y$, and we are considering the quotient of 
$X \ti Y$ that collapses this subspace to the new basepoint, written as $0$ or $[0]$. In other words, 
$ X \me Y $ is the cokernel of the canonical embedding $X \jo Y \to X \ti Y$. The class of the pair $(x, y)$ is 
written as $ x \me y$, and $ x \me 0 = 0 = 0 \me y$.

	The unit of this operation is the discrete space $ \bbS^0 = \{-1, 1\}$, pointed (for instance) at $ 1$.

\vskip -3pt

\LL \begin{small}

	A reader non familiar with this operation might prefer to see now a few computations, in \ref{1.5}.

\end{small} \LB

\vskip 2pt

	A classical counterexample to associativity is based on the triple of spaces $ (\bbQ, \bbQ, \bbN)$, 
proposed by D. Puppe. A proof of the associativity failure, in this case, is given in \cite{MyS}, Section 1.7. 
A class of counterexamples, including the previous one, will be given in Section 4.
	
\Ndt (c) This failure is only concerned with topology. The smash product in $ \Set_\bu$, defined as above, 
is associative up to a canonical bijection, the set-theoretical {\em associator}
    \begin{equation}
\ka\c X \me (Y \me Z) \to (X \me Y) \me Z,   \q   \ka(x \me (y \me z))  =  (x \me y) \me z,
    \label{1.3.2} \end{equation}
so that $ \Set_\bu $ is a symmetric monoidal category, and the forgetful functor 
$ \Top_\bu \to \Set_\bu $ preserves the smash product.

\Ndt (d) We say that a triple $ (X, Y, Z) $ of pointed spaces is {\em associative} (for the smash product) if 
the set-theoretical associator is a (pointed) homeomorphism, that will be called a {\em structural 
homeomorphism}.

	More generally, the terms `structural map' and `structural homeomorphism' of pointed spaces will 
refer to a bijection of the symmetric monoidal closed structure of $\Set\bu$ which happens to be 
continuous or a homeomorphism. A structural homeomorphism will be denoted by the equality sign.

	It is important to note that $\Top\bu$ inherits for free the coherence of those comparisons of 
$\Set\bu$ which happen to be continuous: loosely speaking, they still form commutative diagrams.

	We begin to study the problem, but the main results will be given in the next section.

\subsection{Complements}\label{1.4}
(a) As for the tensor product of modules, the smash product $ X \me Y $ is determined by a universal 
property, the existence of a {\em bipointed} map (i.e.\ pointed in each variable) such that:
    \begin{equation}
\eta\c X \ti Y \todot X \me Y,   \qq   \eta(x, y)  =  x \me y,
    \label{1.4.1} \end{equation}

\vskip -5pt

\LL

\ndt - for every bipointed map $ \ph\c X \ti Y \todot Z $ of pointed spaces there is a unique pointed map 
$ h\c X \me Y \to Z $ such that $ \ph = h\eta$.

\LB

\vskip 3pt

	(Let us note that a bipointed map is always a pointed map, while a bilinear mapping is not 
a homomorphism, generally.)

\Ndt (b) The forgetful functor $ U\c \Top_\bu \to \Top$, which forgets the base point, has a left adjoint
    \begin{equation}
( - )_\bu\c \Top \to \Top_\bu,   \qq   S_\bu  =  S + \sing,
    \label{1.4.2} \end{equation}
which adds to a space an isolated basepoint $*$, by a topological sum. It embeds $\Top$ in $\Top_\bu$ 
as a coreflective subcategory.

\Ndt (c) The functor $( - )_\bu$ transforms the cartesian product of spaces into the smash product of 
pointed spaces (and the unit into the unit)
    \begin{equation}
(S \ti T)_\bu  =  S_\bu \me T_\bu,   \qq   (x, y)  \mapsto  x \me y.
    \label{1.4.3} \end{equation}

\Ndt (d) If the pointed spaces $X, Y$ are compact or Hausdorff, $X \me Y$ is also, as is obvious 
or easily proved. 

\vskip -3pt

\LL

	The smash product of compact Hausdorff pointed spaces is associative, a well-known fact and a 
consequence of \ref{2.2}(b)): their full subcategory has a symmetric monoidal structure. On the other 
hand, locally compact, or exponentiable, or first-countable pointed spaces are not closed under smash 
product, as we shall see in Section 4.

\LB

\subsection{Smash product of compact intervals}\label{1.5}
It is interesting to examine the smash product $X \me Y$ of two non-degenerate euclidean intervals, 
pointed at any point. If $X$ and $Y$ are both compact, the result can be realised in the euclidean plane, 
as is easy to guess and proved below. (The other cases give complex results, that will be examined in 
Section 4.)

	In fact, up to pointed homeomorphism, each of $X$ and $Y$ is either the standard interval $\bbI$ 
or the standard 1-disc $\bbD = \bbD^1 = [- 1, 1]$, both pointed at 0. This gives three spaces $X \me Y$ 
(up to symmetry)
%
%
    
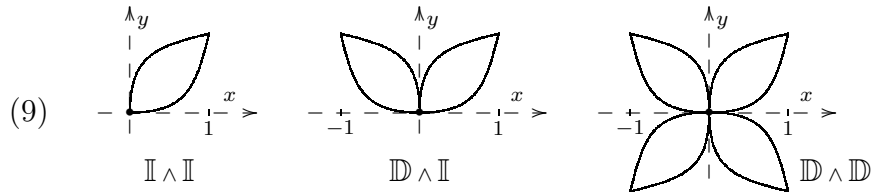
\begin{equation} 
\xy <.35mm, 0mm>:
%
(0,48) *{}; (0,-35) *{};
(-80,-6) *{\sst{1}}; (-72,6) *{\sst{x}}; (-105,34) *{\sst{y}}; 
(-30,-6) *{\sst{-1}}; (30,-6) *{\sst{1}}; (37,6) *{\sst{x}}; (5,34) *{\sst{y}}; 
(80,-6) *{\sst{-1}}; (140,-6) *{\sst{1}}; (147,6) *{\sst{x}}; (115,34) *{\sst{y}}; 
(-95,-22) *{\bbI \me \bbI}; (0,-22) *{\bbD \me \bbI}; (158,-22) *{\bbD \me \bbD}; 
(-110,0) *{\bu}; (0,0) *{\bu};  (110,0) *{\bu};  
(-110,0); (-80,30)  **\crv{(-88,0)&(-85,10)}  **\crv{(-110,22)&(-100,25)},
(0,0); (-30,30) **\crv{(-22,0)&(-25,10)}  **\crv{(0,22)&(-10,25)},
(0,0); (30,30)  **\crv{(22,0)&(25,10)}  **\crv{(0,22)&(10,25)},
(110,0); (80,30) **\crv{(88,0)&(85,10)}  **\crv{(110,22)&(100,25)},
(110,0); (140,30)  **\crv{(132,0)&(135,10)}  **\crv{(110,22)&(120,25)},
(110,0); (80,-30) **\crv{(88,0)&(85,-10)}  **\crv{(110,-22)&(100,-25)},
(110,0); (140,-30)  **\crv{(132,0)&(135,-10)}  **\crv{(110,-22)&(120,-25)},
\POS(-122,0) \ard+(60,0), \POS(-110,-8) \ard+(0,48), 
\POS(-80,-1) \arl+(0,2),  
\POS(-42,0) \ard+(90,0), \POS(0,-8) \ard+(0,48), 
\POS(-30,-1) \arl+(0,2), \POS(30,-1) \arl+(0,2),  
\POS(68,0) \ard+(90,0), \POS(110,-25) \ard+(0,65), 
\POS(80,-1) \arl+(0,2), \POS(140,-1) \arl+(0,2),  
\endxy
    \label{1.5.1} \end{equation}
that can be embedded in the euclidean plane as above, and are respectively homeomorphic to 
the wedge of one, two or four compact discs $\bbD^2$, pointed at a point of the boundary. (Similarly, 
each $n$-ary smash product of non-degenerate compact intervals is a wedge of compact $n$-discs 
$\bbD^n$ pointed as above -- at most $2^n$ of them.)

\vskip 4pt

	The proof, an exercise of topology, can be based on the following map (where $X$ and $Y$ 
are either $\bbI$ or $\bbD$)
    \begin{equation}
f\c X \ti Y \to \bbR^2,   \qq   f(x, y)  =  (|x| \me |y|)(x, y),
    \label{1.5.2} \end{equation}

	The map multiplies the point $(x, y)$ by a scalar factor, its distance from $X \jo Y$. It sends 
$X \jo Y$ to the origin and is injective on its complement, inducing a homeomorphism $X \me Y \to \Im f$ 
of compact Hausdorff spaces.

	In the simplest case $\bbI \me \bbI$, the image of $f$ is the following subspace of the plane
    \begin{equation}
f(\bbI \me \bbI) = \{(x, y) \in X \ti Y  \sep  x^2 \le y \le \sqrt{x}\},
    \label{1.5.3} \end{equation}

\vskip 2pt

\ndt bound by two arcs of parabola, symmetric with respect to the diagonal $y = x$.  The other cases 
follow by symmetry.

\LL \begin{small}

\vskip -2pt

	Indeed, the boundary of $f(\bbI \me \bbI)$ is the union of two arcs, $f(\{1\} \ti \bbI)$ and $f(\bbI \ti \{1\})$. 
The first arc is formed by the points $f(1, t) = (t, t^2)$ of the parabola $y = x^2$ (for  $t \in \bbI$); the second lies 
on $x = y^2$.

\end{small} \LB

\subsection{Symmetry, associativity and smash powers}\label{1.6}
It is not often remarked that the symmetry of the smash product strongly interacts with associativity, 
simplifying this issue.

\Ndt (a) It is easy to see that the triple $ (X, Y, Z) $ is associative if and only if $ (Z, Y, X) $ is.

	In fact, if we assume that the following set-theoretical associator is a homeomorphism
$$
\ka\c X \me (Y \me Z) \to (X \me Y) \me Z,   \;\;\;\;   \ka(x \me (y \me z))  =  (x \me y) \me z,
$$
we get the structural homeomorphism of $ (Z, Y, X) $ as $\psi\ka^{-1}\ph$, using the composed symmetries:
$$
\ph\c Z \me (Y \me X) \to Z \me (X \me Y) \to (X \me Y) \me Z,   \;\;   z \me (y \me x) \mapsto (x \me y) \me z,
$$
$$
\psi\c X \me (Y \me Z) \to (Y \me Z) \me X \to (Z \me Y) \me X,   \;\;   x \me (y \me z) \mapsto (z \me y) \me x.
$$

\Ndt (b) Every triple $ (X, Y, X) $ of pointed spaces is associative: in this case the set-theoretical 
associator is a composed symmetry: 
$$
X \me (Y \me X) \to (Y \me X) \me X \to (X \me Y) \me X.
$$

\Ndt (c) The {\em smash power} $X^{\we n}$ is inductively defined as
    \begin{equation}
X^{\we 0}  =  \bbS^0   \qq   X^{\we n+1}  =  X^{\we n} \we \, X  \qqq   (n \ge 0).
    \label{1.6.1} \end{equation}

	The symmetric formula $X^{\we n+1} = X \we \, X^{\we n}$ can also be used, up to structural 
homeomorphism.

\LL \begin{small}

\vskip -2pt

In fact, assuming that $X^{\we n} \we \, X = X \we \, X^{\we n}$ and applying (b):
$$
X^{\we n+1} \we X = (X^{\we n} \we X) \we X = (X \we X^{\we n}) \we X 
= X \we (X^{\we n} \we X) = X \we X^{\we n+1}.
$$

\end{small} \LB

\Ndt (d) On the other hand, a triple $ (X, X, Y) $ need not be associative: $ (\bbQ, \bbQ, \bbN) $ is not.

\subsection{A colax monoidal structure}\label{1.7}
The smash product of  $\Top_\bu$ can be given a formal status: its finitary extension 
$\, X_1 \me ... \me X_n \,$  has a `colax form' of associativity; more precisely it forms a symmetric 
{\em colax monoidal} category (as is defined in \cite{Ln}, Section 3.1, for the dual lax case).

\LL

\vskip -2pt

	This is briefly investigated in \cite{G2}, Section 5.6. The extended operation can be defined in 
any pointed category with finite limits and colimits: see  \cite{CrJ}, Definition 4.1.

\LB

\skp

	The (regular) $n$-{\em ary smash product} $\, X_1 \me ... \me X_n \,$ of pointed spaces is defined 
by the universal multi-pointed map (i.e.\ pointed in each variable)
    \begin{equation}
\eta\c X_1 \ti ... \ti X_n \; \todot X_1 \me ... \me X_n,   \q   \eta(x_1, ..., x_n)  =  x_1 \me ... \me x_n.
    \label{1.7.1} \end{equation}

	The solution is a quotient map, the cokernel of the embedding of the subspace $H$ of coordinate 
hyperplanes
    \begin{equation}
H \to X_1 \ti ... \ti X_n,   \q   H  =  \{(x_i) \in \Pro X_i \sep  \exists i\c  x_i = 0\}.
    \label{1.7.2} \end{equation}

	The term $ x_1 \me ... \me x_n $ annihilates if and only if some coordinate $ x_i $ is zero; all the 
other equivalence classes are singletons.

	As described in \cite{Ln} (in the dual lax case), there is a system of {\em associativity comparisons}, 
as in the following example
    \begin{equation}
\ga\c X_1 \me ... \me X_5 \; \to \; (X_1 \me (X_2 \me X_3 \me X_4)) \me X_5.
    \label{1.7.3} \end{equation}

\LL \begin{small}

\vskip -2pt

	Each comparison is produced by the universal property of $\eta$, and determined by an 
arrangement of parentheses, more formally a tree.

\end{small} \LB

\skp
	
	Here the {\em coherence} of the comparisons (loosely speaking, the fact that they produce 
commutative diagrams) is automatically granted by the underlying monoidal structure of $\Set_\bu$, 
as in \ref{1.3}(d).
	
	Of course, these comparisons need not be invertible. (In the general colax case, if all of them are 
we go back to a monoidal structure in {\em unbiased} form, that is based on finite tensor products.)

\subsection{Regular associativity}\label{1.8}
(a) In particular, a triple of pointed spaces $(X, Y, Z)$ has two comparisons
    \begin{equation} \begin{array}{cl}
\ga'\c X \me Y \me Z \to X \me (Y \me Z),   &\;\;   \ga'(x \me y \me z)  =  x \me (y \me z),
\\[3pt]
\ga''\c X \me Y \me Z \to (X \me Y) \me Z,   &\;\;  \ga''(x \me y \me z)  =  (x \me y) \me z,
    \label{1.8.1} \end{array} \end{equation}
produced by the universal property of the tri-pointed mapping $ \eta\c X \ti Y \ti Z \todot  X \me Y \me Z$.

	Composing with the latter we get two surjective tri-pointed maps $\rho' = \ga'\eta$ and 
$\rho'' = \ga''\eta$, defined on the cartesian product $X \ti Y \ti Z$
    \begin{equation} \begin{array}{c} 
    \xymatrix  @C=10pt @R=20pt
{
X \ti Y \ti Z~   \ar[r]^-{\eta}   \ar[d]_(.45){1 \times p'}   \ar[rd]|-{~\rho'~}    &   
~X \me Y \me Z~    \ar[d]^(.45){\ga'}  &
X \ti Y \ti Z~   \ar[r]^-{\eta}   \ar[d]_(.45){1 \times p''}   \ar[rd]|-(.6){\rho''}    & 
~X \me Y \me Z   \ar[d]^(.45){\ga''} 
\\ 
X \ti (Y \me Z)~   \ar[r]_-{q'}    &   ~X \me (Y \me Z)~  &
(X \ti Y) \me Z~   \ar[r]_-{q'}    &   ~(X \me Y) \me Z
}
    \label{1.8.2} \end{array} \end{equation}
    \begin{equation*}
\;\;\;\;   \rho'(x, y, z)  =  x \me (y \me z),     \qq    \rho''(x, y, z)  =  (x \me y) \me z.
    \label{1.8.2bis} \end{equation*}

\Ndt (b) With some abuse of terminology, we say that a particular instance $X \me (Y \me Z)$ is a 
{\em regular smash product} if its comparison $\ga'$ is invertible (a property that makes sense 
in a general lax or colax monoidal category). Each of the following conditions is equivalent to the 
former

\LL

\ndt  -  $\rho'\c X \ti Y \ti Z  \to  X \me (Y \me Z)$ is a quotient map,

\Ndt  - $1 \ti p'$ sends every open (or closed) subset of $X \ti Y \ti Z$ {\em saturated for} $\rho'$ 
to an open (or closed) subset of $X \ti (Y \me Z)$.

\LB

\skp

	Similarly we say that $(X \me Y) \me Z$ is a regular smash product if the same properties hold 
for the maps $\ga''$, $\rho''$ and $1 \ti p''$, in the right diagram \eqref{1.8.2}.

\LL \begin{small}
\vskip -2pt

	One could say that the triple $(X, Y, Z)$ is {\em left regular} or {\em right regular}, which would 
be formally correct but liable to confusion.

\vskip -2pt
\end{small} \LB

\skp

	In the first case $\rho'$ is a quotient map, and the canonical bijection $\ka$ (of sets) of 
the following diagram is continuous, hence a structural map in $\Top_\bu$ 
    \begin{equation} \begin{array}{ccc} 
    \xymatrix  @C=15pt @R=0pt
{
~X \ti Y \ti Z~   \ar[r]^-{\rho'}   \ar[rdd]_-{\rho''}    &   ~X \me (Y \me Z)~     \ard[dd]^-{\ka} 
\\ 
&&    \ka(x \me (y \me z))  =  (x \me y) \me z.
\\ 
&   ~(X \me Y) \me Z~ 
}
    \label{1.8.3} \end{array} \end{equation}

	In the second case there is a structural map $\ka'$ the other way round.

\Ndt (c) If both occurrences are regular, $\ka$ is a structural homeomorphism and we say that the triple 
$ (X, Y, Z) $ is {\em regularly associative}. This is strictly stronger than associativity.

\LL \begin{small}

\vskip -6pt

For instance, the triple $ (\bbQ, \bbN, \bbQ)$ is associative, by symmetry, but is not regularly 
associative:  $ (\bbQ \me \bbN) \me \bbQ = \bbQ \me (\bbQ \me \bbN)$ is not a regular 
smash product (see Theorem \ref{4.8}).

\end{small} \LB

\vskip 4pt

	It is also useful to note that, if one of $ X \me (Y \me Z) $ and $ (X \me Y) \me Z $ is a regular 
smash product and the other is not, the triple $(X, Y, Z)$ is not associative.

\Ndt (d) Concretely, the ternary smash products $X \me (Y \me Z)$ and $(X \me Y) \me Z$ 
can be viewed as two topologies on the quotient {\em set}
$$
(X \ti Y \ti Z) / (X \jo Y \jo Z),
$$
each of them (weakly) coarser than the quotient topology, and we are considering whether one of 
them coincides with the latter, or both.

\Ndt (e) {\em Note}. We already mentioned that the article \cite{CrJ} studies the smash product in a 
pointed category, say $\mathsf{A}$, with finite limits and colimits; the property of 
$\me$-{\em associativity} of $\mathsf{A}$ is introduced in Definition 4.2. In the present terminology, 
extended to the general situation, this means that all triples are regularly associative.

	The weaker form of associativity of \ref{1.3}(d) can also be extended if the forgetful functor 
$\mathsf{A} \to \Set_\bu$ `represented' by the unit of the smash product is faithful.

\section{Exponentiable pointed spaces and function spaces}\label{s2}
	Exponentiable pointed spaces and function spaces in $ \Top_\bu $ are defined with respect to the 
smash product. (We recall that, in a pointed category, the only exponentiable object for the cartesian 
product is the zero object.)

	Exponentiability is closely related to partial associativity properties of the smash product, as shown 
by Theorem \ref{2.2}.

\subsection{Exponentiable pointed spaces}\label{2.1}
(a) The pointed space $ A $ is said to be {\em exponentiable} if the functor 
$- \me A\c \Top_\bu \to \Top_\bu$ has a right adjoint. The latter, written as $ (-)^A$, is characterised by 
a natural bijection $\ph_{XY}$, the exponential law of pointed spaces
    \begin{equation}
- \me A \, \adj \, ( - )^A,     \q     \ph_{XY}\c \Top_\bu(X \me A, Y) \to \Top_\bu(X, Y^A).
    \label{2.1.1} \end{equation}

\Ndt (b) If the underlying space $ |A| $ is exponentiable (for the cartesian product), then it is easy to see 
that $A$ is exponentiable in $\Top_\bu$ (for the smash product): the exponential $ Y^A $ is the set of 
pointed maps $ \Top_\bu(A, Y) $ with the topology induced by the exponential topology of the whole 
space of unpointed maps $ |Y|^{|A|}$. The natural bijection $\ph_{XY}$  (corresponding to the structure 
of $\Set_\bu$) is induced by  the unpointed one, in $ \Top $ (recalled in \eqref{1.1.1}).

	It is also true that a pointed space which is exponentiable in $ \Top_\bu $ for the smash product 
has an underlying space exponentiable in $\Top$, as proved by Cagliari \cite{Ca}. The terms 
{\em exponentiable pointed space} and {\em pointed exponentiable space} are thus equivalent and 
non-ambiguous.

\Ndt (c) Let $A$ and $B$ be exponentiable pointed spaces. Then $A \ti B$ is also, by \ref{1.1}(c), but 
$A \me B$ `rarely' is, unless $A$ and $B$ are both compact: for instance, $\bbR \me \bbR$ and 
$\bbR \me \bbI$ are not (see Theorem \ref{4.2}).

	However, it is still true that the endofunctors $( - )^A$ and $( - )^B$ commute: they are subfunctors 
of the corresponding unpointed exponentials, which commute by \ref{1.1}(c). Their left adjoints  
$- \me A$ and $- \me B$ also commute (using symmetries, this follows from \ref{2.2}(a)).

\Ndt (d) We now study the relationship of associative triples with exponentiable or compact 
exponentiable pointed spaces. The main results are gathered in the following theorem, which partially 
depends on a property of compactness deferred to \ref{2.4}.

\subsection{Theorem {\rm (Exponentiability and regular associativity)}}\label{2.2}
{\em
Let $X$, $Y$, $Z$ be pointed spaces.

\Ndt (a) If $ X $ and $ Z $ are exponentiable (e.g.\ locally compact), the triple $ (X, Y, Z) $ is regularly 
associative. In particular, the set-theoretical associator $ \ka\c X \me (Y \me Z) \to (X \me Y) \me Z $ is 
a structural homeomorphism.

\Ndt (b) If $ Y $ is compact and $ X $ or $ Z $ is compact and exponentiable (e.g.\ compact Hausdorff), 
the triple $ (X, Y, Z) $ is regularly associative.

\Ndt (c) If $ Y $ and $ Z $ are compact and exponentiable, $ Y \me Z $ is also and there are structural 
homeomorphisms (given by the canonical bijections of $\Set_\bu$), for all pointed spaces $X$:

\vskip -6pt

    \begin{equation}
(X \me Y) \me Z  =  X \me (Y \me Z)  =  X \me (Z \me Y)  =  (X \me Z) \me Y,
    \label{2.2.1} \end{equation}
    \begin{equation}
(X^Z)^Y  =  X^{Y \me Z} =  X^{Z \me Y} =  (X^Y)^Z.
    \label{2.2.2} \end{equation}

\vskip 5pt

	In particular,  the exponentials $( - )^Y$ and $( - )^Z$ commute. One can replace $Y$ and $Z$ 
with their smash powers $Y^{\we m}$ and $Z^{\we n}$, which are also compact and exponentiable. 
}
\begin{proof}
(a) If $ X $ is exponentiable, the functor $ X \ti - $ preserves the quotient $p'\c Y \ti Z \to Y \me Z$
in \eqref{1.8.2} and $ X \me (Y \me Z) $ is a regular smash product. If also $ Z $ is exponentiable, 
$ (X \me Y) \me Z $ is regular as well.

\Ndt (b) First we suppose that $Y$ and $Z$ are compact. Applying property \ref{2.4} to the pointed space 
$ H = Y \ti\, Z $ and its compact subspace $ Y \jo Z$, the map $ 1 \ti p\c X \ti Y \ti Z \to X \times (Y \me Z) $ 
is a quotient map, and $ X \me (Y \me Z) $ is a regular smash product. 

	Adding the hypothesis that $ Z $ is exponentiable, $ (X \me Y) \me Z $ is regular as well, as in the 
previous point. The other case is symmetrical.

\Ndt (c) Now $ Y $ and $ Z $ are supposed to be compact and exponentiable. From (b) and symmetry 
we get \eqref{2.2.1}. As to \eqref{2.2.2}, the composed endofunctor $ F = (- \me Y) \me Z $ has a 
composed right adjoint $ G = (( - )^Z)^Y$. But $ F = - \me (Y \me Z)$, by \eqref{2.2.1}, which implies 
that $ Y \me Z $ is exponentiable and $ G = ( - )^{Y \me Z}$.

	The last point is also a consequence of (b).
\end{proof}
%

\subsection{Lemma {\rm (Exponentiability and associativity)}}\label{2.3}
{\em
If $ Y, Z $ are exponentiable pointed spaces and all triples $ (X, Y, Z) $ are associative, the pointed space 
$ Y \me Z $ is exponentiable.
}
\begin{proof}
The composed endofunctor $ F = (- \me Y) \me Z $ has a composed right adjoint $ G = (( - )^Z)^Y$. 
Now $ F = - \me (Y \me Z) $ implies that $ Y \me Z $ is exponentiable. Moreover $ G = ( - )^{Y \me Z}$, 
but this is already in \ref{2.2}(c).
\end{proof}
%

\subsection{A property of compactness}\label{2.4}

We have already used the following fact: if $K$ is a compact non-empty subspace of the space $H$, 
the surjective map $f = 1 \ti p\c X \ti H \to X \ti (H/K)$ is a quotient map, for every space $X$.

	This is Theorem 6.2.4(b) of Maunder's book \cite{Mn}, a textbook not frequently available. 
We write down a (slightly shorter) proof.

	Let $ W \sub X \times H $ be an open $f$-saturated subset. We have to prove that $ f(W) $ is open in 
$ X \ti (H/K)$. Equivalently, every $ (x, y) \in W $ has a basic nbd $ U \ti V \sub W $ (a product of open 
sets) such that $ f(U \ti V) = U \ti p(V) $ is open in $ X \ti (H/K)$; this will be certainly true if $ V $ is 
$p$-saturated.

	There are two cases.

\Ndt (a) If $ y \in K$, then $ \{x\} \ti K \sub W$, because $ W $ is $f$-saturated. The subspace $ \{x\} \ti K $ 
is covered by a family of basic nbds $ U_z \ti V_z \sub W $ of the points $ (x, z) \in \{x\} \ti K$, and there is 
a finite subfamily $ U_i \ti V_i \sub W$, where $ (V_i) $ still covers $ K$. We have thus an open nbd 
$ U = \Cap_i \, U_i $ of $ x $ and an open subset $ V = \Cup_i \, V_i \sups K$. Finally $ U \ti V $ is an 
open nbd of $ (x, y) $ contained in $W$
$$
U \ti V   \; = \;   \Cup_i \, U \ti V_i  \; \sub \;   \Cup_i \, U_i \ti V_i   \; \sub \;   W,
$$
and $ V $ is $p$-saturated, as it contains $ K$.

\Ndt (b) If $ y \notin K$, there is an open nbd $ U \ti V $ of $ (x, y) $ in $ W$. If $ V $ does not meet $ K$, it is 
$p$-saturated. Otherwise, let $ y' \in V \cap K$. We apply the previous point to $ (x, y') \in U \ti V \sub W$, 
and we find an open nbd $ U' \ti V' \sub W $ with $ V' \sups K$. Then $ (U \cap U') \ti (V \cup V') $ is an 
open nbd of $ (x, y) $ contained in $ W$, with $ V \cup V' \sups K$

\section{Homotopy constructions}\label{s3}
	The cylinder functor, path functor, cones and suspension for pointed spaces can be defined 
adapting the unpointed case.

	But cylinder, cones and suspension can be equivalently obtained in an `internal way', using the 
smash product (with fixed well-behaved spaces). All these homotopy constructions have a right adjoint: 
path space, cocones and loop spaces, all of them obtained as pointed exponentials.

	Part of these facts are well known in topology, but their interaction with the smash product is 
often omitted, or confined to convenient spaces. For a categorist, this section is meant to show part of 
the interest of the smash product.

\subsection{Pointed cylinder and path space}\label{3.1}
Homotopies of pointed spaces are ordinary homotopies of pointed maps that keep the basepoint fixed. 
They are represented by the {\em cylinder functor of pointed spaces}, or {\em pointed cylinder}, or 
{\em reduced cylinder}.

	On the pointed space $ X $ this is defined as the quotient $ IX $ of the unpointed cylinder 
$ |X| \ti \bbI $ that collapses the fibre at $ 0_X$, forming the new basepoint
    \begin{equation}
IX  =  (|X| \ti \bbI)/(\{0_X\} \ti \bbI),   \q   0_{IX}  =  [0_X, t]   \q\;\;   (t \in \bbI),
    \label{3.1.1} \end{equation}
so that a (pointed) map $ IX \to Y $ amounts to a continuous mapping $ |X| \ti \bbI \to |Y| $ sending all 
pairs $ (0_X, t) $ to $ 0_Y$.

\skp	More simply, the {\em path space} of a pointed space 
    \begin{equation}
P(Y)  =  (|Y|^\bbI, \om_0),   \q   \om_0(t)  =  0_Y   \q\;\;   (t \in \bbI),
    \label{3.1.2} \end{equation}
is the path space $ |Y|^\bbI $ of the underlying topological space (with compact-open topology), 
pointed at the trivial loop $ \om_0 $ of the basepoint.

	We have thus the {\em path functor} $P$ of pointed spaces, right adjoint to the pointed 
cylinder ($I \adj P$)
    \begin{equation}
P\c \Top_\bu \to \Top_\bu,  \q   P(g)(a)  =  ga   \q\;\;   (g\c Y \to Y', \;  a \in |Y|^\bbI).
    \label{3.1.3} \end{equation}

	The adjunction says that, in $ \Top_\bu$, a homotopy of maps $ X \to Y $ can be equivalently 
described by a map $\hat{\ph}\c IX \to Y $ or a map $\check{\ph}\c X \to PY$.

\subsection{The pointed interval}\label{3.2}
The {\em pointed interval} comes out of the standard interval $ \bbI$, adding a basepoint as in 
\ref{1.4}(b)
    \begin{equation}
\bbI_\bu =  \bbI + \sing,
    \label{3.2.1} \end{equation}
\ndt  a compact Hausdorff pointed space. The cylinder and path functor of $ \Top_\bu $ can thus be 
expressed as the adjoint functors
    \begin{equation} \begin{array}{ccc}
I\c \Top_\bu \to \Top_\bu,   &\q\;\;&   IX  =  X \me \bbI_\bu,
\\[5pt]
P\c \Top_\bu \to \Top_\bu,   &&   PY  =  Y^{\bbI_\bu},
    \label{3.2.2} \end{array} \end{equation}
using the structure of $ \Top_\bu$, smash product and pointed exponentials, instead of the cartesian 
structure of $ \Top$, as previously.

	Applying \ref{1.4}(c), the smash-powers of $\bbI_\bu$ correspond to the standard cubes $\bbI^n$
    \begin{equation}
(\bbI_\bu)^{\we n}  =  (\bbI^n)_\bu,
    \label{3.2.3} \end{equation}
and the associativity property of \ref{2.2}(c) gives:
    \begin{equation} \begin{array}{c}
I^nX  =  (...(X \me  \bbI_\bu) ...) \me  \bbI_\bu  =  X \me (\bbI_\bu)^{\we n}
\\[5pt]
\qq\;\;   =  X \me (\bbI^n)_\bu  =  (|X| \times \bbI^n) / (\{0_X\} \times \bbI^n).
    \label{3.2.4} \end{array} \end{equation}

	Moreover, by \ref{2.2}(c), the iterated path space 
    \begin{equation}
P^n Y  =  Y^{{(\bbI_\bu)^{\we n}}}  =  Y^{(\bbI^{\we n})_\bu},
    \label{3.2.5} \end{equation}
is the set $ \Top_\bu((\bbI^n)_\bu, Y) = \Top(\bbI^n, Y) $ with the compact-open topology, pointed at 
the constant map at $ 0_Y$.

\subsection{Pointed cones}\label{3.3}
The {\em upper cone} $ C^+X $ of a pointed space is the cokernel of the upper face $ \ddp $ of the 
pointed cylinder $ IX$
    \begin{equation}
C^+X  =  \Coker (\ddp\c X \to IX),   \q\;\;   \ddp(x)  =  [x, 1].
    \label{3.3.1} \end{equation}

	The canonical projection $ \ga\c IX \to C^+X $ collapses the upper basis $\ddp X$ to an upper 
vertex $ v^+$. This gives a functor $ C^+ $ and a natural transformation
    \begin{equation}
C^+\c \Top_\bu \to \Top_\bu,   \q\;\;   \ga\c I \to C^+,
    \label{3.3.2} \end{equation}
characterised by the following universal property:

\LL

\vskip -2pt
\ndt - for every (pointed) homotopy $ \ph\c f \simeq 0 $ reaching the zero-map $ X \to Y$, there is precisely 
one map $ h\c C^+X \to Y $ such that $ \ph = h\ga\c IX \to C^+X \to Y$.

\LB

\skp

	Dually, the {\em lower cone} $ C^-X = \Coker (\ddm\c X \to IX) $ is obtained by collapsing the 
lower basis of $ IX $ to a lower vertex $ v^-$.

	The cones of the pointed space $ X $ can be obtained as smash products
    \begin{equation}
C^\al X  =  X \me \bbI_\al  =  (X \ti \bbI) / (X \ti \{\al\} \cup \{0\} \ti \bbI)   \q   (\al = 0, 1),
    \label{3.3.3} \end{equation}
where $ \bbI_\al $ is the standard interval, {\em pointed at} $\al$.

	Although the upper and lower cones are obviously homeomorphic, it is often convenient to 
distinguish them. In Directed Algebraic Topology this is necessary \cite{G1}.

\subsection{Pointed suspension}\label{3.4}
The (pointed) {\em suspension} $ \Si X $ of a pointed space, also called {\em reduced suspension}, 
can be obtained as the colimit of the (solid) left diagram below, using the pointed cylinder $ IX $ 
and its faces $ \dd^\al(x) = (x, \al)$
%
    \begin{equation}
    \begin{array}{c}  \xymatrix  @C=11pt @R=14pt
{
&  ~X~    \ar[r]  \ar[d]_-{\ddp}      &  \sing    \ard[dd]^-{v}  &&
&  ~X~    \ar[r]  \ar[d]_-{\ddp}      &  \sing    \ard[d]^-{v}  
\\  
~~X~~   \ar[r]^-{\ddm}    \ar[d]  &  ~IX~   \ard[rd]^-{\si}  &&&
~~X~~   \ar[r]^-{\ddm}    \ar[d]  &  ~IX~   \ard[r]  \ard[d]    &  ~C^+(X)~   \ard[d]
\\  
~~\sing~~   \ard[rr]_-{v}    &&  ~\Si X~   &&
~~\sing~~   \ard[r]    &   ~~C^-(X)~~   \ard[r]  &  ~\Si X~ 
}
    \label{3.4.1} \end{array} \end{equation}

	Equivalently it is produced by three pushouts, in the right diagram. The map $ v\c \sing \to \Si X $ 
is the zero map. (The vertices $ v^-, v^+ $ of the unpointed suspension are here identified.)

	Pasting the lower pushouts (resp.\ the right-hand pushouts) in the second diagram above, 
$ \Si X $ is bound to the cones $ C^\al(X) $ by the following pushouts
%
    \begin{equation}
    \begin{array}{c}  \xymatrix  @C=5pt @R=7pt
{
~~X~~   \ar[rr]^-{\ddm}  \ar[dd]   &&   ~C^+(X)~   \ar[dd]   &&&
~~X~~   \ar[rr]  \ar[dd]_-{\ddp}   &&  ~~\sing~~   \ar[dd]
\\
&&     \arld@/_/[ld]   &&&&&     \arld@/_/[ld]
\\
~~\sing~~   \ar[rr]  &&  ~\Si X~     &&&   ~C^-(X)~   \ar[rr]  &&  ~\Si X~ 
}
    \label{3.4.2} \end{array} \end{equation}

\Ndt showing that $ \Si X $ is the quotient of each cone that collapses its basis to a point
    \begin{equation} \begin{array}{c}
\Si X  =  \Coker (\ddm\c X \to C^+(X))  =  C^+(X) / \ddm X
\\[5pt]
\q   =  \Coker (\ddp\c X \to C^-(X))  =  C^-(X) / \ddp X.
    \label{3.4.3} \end{array} \end{equation}

	We have thus the pointed suspension functor $\Si$ and a natural transformation
    \begin{equation}
\Si \c \Top_\bu \to \Top_\bu,   \qq   \si\c I \to \Si,
    \label{3.4.4} \end{equation}
characterised by the following universal property:

\LL

\vskip -2pt

\ndt - for every (pointed) endohomotopy $ \ph\c 0 \eq 0 $ of the zero map $ X \to Y $ there is precisely 
one map $ h\c \Si X \to Y $ such that $ \ph = h\si\c IX \to \Si X \to Y$.

\LB

\Ndt {\em Remarks}. (a) The suspension of pointed spaces can also be expressed as a smash product
    \begin{equation}
\Si X  =  X \me \bbS^1.
    \label{3.4.5} \end{equation}

\Ndt (b) The following endofunctors of $ \Top_\bu$
$$
I = - \me \bbI_\bu,   \q   C^+ = - \me \bbI_+,   \q   C^- = - \me \bbI_-,   \q   \Si  = - \me \bbS^1,
$$
commute with each other, and with each functor $- \me A$ produced by an exponentiable 
pointed space, by \ref{2.1}(c).

\Ndt (c) The suspension $ \Si X $ can also be obtained pasting two cones of $ X $ along their faces.

\subsection{The pointed spheres. }\label{3.5}
The pointed $n$-sphere can be defined as the topological quotient of the euclidean $n$-cube $\bbI^n$
    \begin{equation}
\bbS^n  =  \bbI^n/\dd \bbI^n   \qq   (n > 0),
    \label{3.5.1} \end{equation}
that collapses the boundary of $ \bbI^n $ in $ \bbR^n $ to the basepoint.

	The pointed $n$-sphere is a smash power of the pointed circle
    \begin{equation}
\bbS^n  =  \bbS^1 \me ...  \me \bbS^1  =  (\bbS^1)^{\we n},
    \label{3.5.2} \end{equation}
as a consequence of Remark (a), below.

	Applying \ref{2.2}(c), the iterated pointed suspension $ \Si ^nX $ is produced by the smash 
product with $ \bbS^n $
    \begin{equation}
\Si ^nX  =  (...(X \me \bbS^1)...) \me \bbS^1  =  X \me (\bbS^1)^{\we n} =  X \me \bbS^n,
    \label{3.5.3} \end{equation}
and all spheres are suspensions of $\bbS^0$
    \begin{equation}
\Si ^n\bbS^0  =  \bbS^0 \me \bbS^n  =  \bbS^n.
    \label{3.5.4} \end{equation}
{\em Remark.} (a) For arbitrary spaces $ X, Y $ with non-empty subspaces $ H, K$
    \begin{equation}
X/H \me Y/K  =  (X \ti Y)/(X \ti K \cup H \ti Y),
    \label{3.5.5} \end{equation}
because the composed projection $X \ti Y \to X/H \times Y/K \to X/H \me Y/K$ collapses the 
subspace $X \ti K \cup H \ti Y$ to the basepoint.

\LL \begin{small}

\vskip -2pt

	This also works with empty subspaces, provided that 
the pointed space $ X/H $ is defined by the left adjoint to the obvious inclusion 
$ \Top_\bu \to \Top_2 $ in the category of relative pairs of topological spaces, so that $ X/\es = X_\bu$.

\end{small} \LB

\subsection{Pointed cocones and loop space}\label{3.6}
For a pointed space $ X$, the cocones $ E^\al X $ and the loop-space $\Om X $ are defined by 
pullbacks, or kernels, in $ \Top_\bu$ (with $\al = 0, 1$)
%
    \begin{equation}
    \begin{array}{c}  \xymatrix  @C=3pt @R=7pt
{
~E^\al(X)~   \ar[rr]  \ar[dd]   &&   ~~\sing~~   \ar[dd]   &&
~\Om X~   \ar[rr]  \ar[dd]   &&  ~~\sing~~   \ar[dd]   &&
~\Om X~   \ar[rr]  \ar[dd]   &&  ~~\sing~~   \ar[dd]
\\
\arld@/_/[ru]   &&&&     \arld@/_/[ru]   &&&&     \arld@/_/[ru]
\\
~P(X)~~   \ar[rr]_-{\dd^\al}  &&  ~~X~~     &&
~E^+(X)~   \ar[rr]_-{\ddm}  &&  ~~X~~     &&
~E^-(X)~   \ar[rr]_-{\ddp}  &&  ~~X~~ 
}
    \label{3.6.1} \end{array} \end{equation}
    \begin{equation} \begin{array}{l}
E^\al X  =  \{a \in PX  \sep  a(\al) = 0\}  =  \Ker (\dd ^\al\c PX \to X),
\\[5pt]
\, \Om X  =  \{a \in PX  \sep  a(0) = 0 = a(1)\}
\\[2pt]
\q\,    =  \Ker (\ddm\c E^+X \to X)  =  \Ker (\ddp\c E^-X \to X).
    \label{3.6.2} \end{array} \end{equation}

	The following endofunctors of $ \Top_\bu$, right adjoints to those of \ref{3.4}(b)
    \begin{equation}
P  =  (- )^{\bbI_\bu},   \q   E^+  =  ( - )^{\bbI_+},   \q   E^-  =  ( - )^{\bbI_-},   \q   \Om  =  ( - )^{\bbS^1},
    \label{3.6.3} \end{equation}
commute with each other, and with each functor $ ( - )^A$ produced by an exponentiable 
pointed space.

\LL \begin{small}

\vskip -2pt

	We recall that in $ \Top $ cones and suspension do not preserve sums and do not have a 
right adjoint.

\end{small} \LB

\section{Non-associative triples}\label{s4}
	The associativity formula $ X \me (Y \me Z) = (X \me Y) \me Z $ certainly holds when $ X $ 
and $ Z $ are exponentiable. But it can fail when $ X $ and $ Y $ (or $ Y $ and $ Z$) are so.

	For instance, it is the case for any pair $X, Y$ of non-degenerate pointed intervals where $X$ is 
not compact. This case is studied in the first part of this section.

	The second part examines a variety of non-associative triples, generalising the previous cases and 
also the usual example $(\bbQ, \bbQ, \bbN)$.

\subsection{Smash product of euclidean intervals}\label{4.1}
Examining the smash product $X \me Y$ of two non-degenerate euclidean intervals, we already 
considered the compact case in \ref{1.5}. We suppose now that at least one of these intervals is not 
compact.

	Then $X \me Y$ is still Hausdorff, but it is neither locally compact (by Theorem \ref{4.2}) 
nor first countable at the basepoint (by Proposition \ref{4.3}). It cannot be realised in the plane, 
nor in any metric space, but is nevertheless an interesting space which can be accurately described.

	As a consequence, $X \me Y$ is not exponentiable \cite{He}, and there exists some pointed 
space $Z$ such that:

\LL

\vskip -5pt

\ndt - the associativity formula $ X \me (Y \me Z) = (X \me Y) \me Z$ fails (by \ref{2.3}),

\vskip 2pt

\ndt - the smash product $(X \me Y) \me Z$ fails to be regular.

\LB

\vskip 4pt

	(The second point follows from the first: $X$ is exponentiable and the smash product 
$ X \me (Y \me Z) $ is regular.) 

	More explicitly, Theorem \ref{4.8} shows that these properties fail for $Z = \bbQ$, the rational line. 
(The proof relies on a sequence $(r_n)$ of rational numbers that converges to an irrational number, 
taking into account that all sequences $(r_n/m)_{n\in \bbN}$ also do, for $m > 0$.)

\subsection{Theorem}\label{4.2}
{\em
Let $X$ and $Y$ be a pair of non-degenerate pointed intervals. The smash product $X \me Y$ is 
locally compact (and exponentiable) if and only if $X$ and $Y$ are both compact.

	Otherwise the basepoint has no compact nbd; at any other point $X \me Y$ is locally euclidean.
}
\begin{proof}
The compact case is already known, from \ref{1.5}.

	We write as $ p\c X \ti Y \to X \me Y $ the canonical projection. We can assume that both intervals 
are pointed at 0. We also assume that $X$ is non-compact, which gives three cases (up to pointed 
homeomorphism), that is $[0, 1[$, $[-1, 1[$  and $]-1, 1[$; in any case $1 \notin X \supset [0, 1[$.

\vskip 4pt

	The following figures represent the space $X \ti Y$  when 0 is an endpoint of both intervals
%
%
    \begin{equation} 
\xy <.4mm, 0mm>:
%
(0,25) *{}; (0,-22) *{};
(36,1) *{\sst{C}}; (70,-12) *{[0, 1[^2}; (30,-5) *{\bu}; (37,-9) *{\bu}; (43,-13) *{\bu}; 
(136,1) *{\sst{C}}; (175,-12) *{[0, 1[ \times \bbI}; (130,-5) *{\bu}; (137,-9) *{\bu}; (143,-13) *{\bu}; 
\POS(10,-20) \arl+(39,0) \arl+(0,39), \POS(50,20) \arld-(40,0) \arld-(0,40),
\POS(10,-19.6) \arl+(39,0),  \POS(10.4,-20) \arl+(0,39);
\POS(110,-20) \arl+(39,0) \arl+(0,39), \POS(150,20) \arl-(40,0) \arld-(0,40),
\POS(110,-19.6) \arl+(39,0),  \POS(110.4,-20) \arl+(0,39);
\endxy
    \label{4.2.1} \end{equation}

\LL \begin{small}

	The thick lines represent the subspace $X \jo Y$ which is collapsed to the basepoint  [0]  in
 the quotient $X \me Y$.
	
\end{small} \LB

\skp

	Taking on the general case, an arbitrary nbd of $ [0] $ in $ X \me Y $ is the image $ p(W) $ of a 
nbd $ W $ of $ X \jo Y $ in $ X \ti Y$, and we want to show that $ p(W) $ cannot be compact, constructing 
a closed non-compact subset $ p(C) $ as in the figure above.
	
	We fix a strictly increasing sequence $ (x_n) $ in $ X $ that converges to $1 $ in $ [1/2, 1]$, like
    \begin{equation}
x_n  =  (n +1)/(n +2)   \qq   (n \ge 0).
    \label{4.2.2} \end{equation}

	It will be repeatedly used in the next subsections.

\vskip 2pt

	For every $ n \ge 0$, $(x_n, 0) \in X \ti \{0\}$, and we choose a point 
$ w_n = (x_n, y_n) \in W \setm (X \jo Y)$, with $ 0 < |y_n| < 1/(n +1)$. The support 
$ C = \{w_n \sep n \ge 0\} $ of this sequence is closed in $ X \ti Y $ (it is the trace of its closure 
$\overline{C} = C \cup \{(1, 0)\} $ in $ \bbI^2$). 

	Now the discrete space $ p(C) $ is closed in $X \me Y$ (because $ C $ does not meet $ X \jo Y$) 
and in $ p(W)$. The compactness of the latter would imply the compactness of $p(C)$, which is discrete 
and infinite.

	Finally, the complement of the basepoint of $X \me Y$ is homeomorphic to an open subspace of 
the plane, namely $(X \setm \{0\}) \times (Y \setm \{0\})$, by a restriction of $p$.
\end{proof}
%

\subsection{Proposition}\label{4.3} 
{\em
Let $X$ and $Y$ be a pair of non-degenerate pointed intervals. The smash product $X \me Y$ is first 
countable if and only if $X$ and $Y$ are both compact.

	Otherwise $X \me Y$ has no countable local basis at the basepoint, and cannot be embedded in 
any first-countable space.

	In particular $[0, 1[ \me \bbI$  cannot be embedded in $\bbI \me \bbI$,  showing that the functor 
$- \me \bbI$ does not preserve subspaces, generally.
}
\begin{proof}
As in the proof of Theorem \ref{4.2}, the compact case is already known and we can assume that  
$1\notin X \supset [0, 1[$. Given a sequence of open nbds $ (W_n) $ of $ X \jo Y $ in $X \ti Y$, we want
to construct an open nbd $ W $ that does not contain any of them.

	We begin by forming an open nbd $ W'_\ep $ of $ X \ti \{0\} $ in $ X \ti Y$, as in the following picture
%
%
    \begin{equation} 
\xy <.4mm, 0mm>:
%
(-15,27) *{}; (0,-30) *{};
(75,-25) *{\sst{x_0}}; (78,15) *{\sst{w_0}}; (75,10) *{\bu}; 
(100,-25) *{\sst{x_1}}; (103,5) *{\sst{w_1}}; (100,0) *{\bu}; 
(-7,10) *{\sst{\ep_0}}; (-7,0) *{\sst{\ep_1}}; 
(170,-15) *{X \ti Y}; (40,-5) *{W'_\ep};
\POS(0,-20) \arl+(149,0) \arl+(0,39), \POS(150,-20) \arld+(0,40), 
\POS(0,-19.6) \arl+(149,0),  \POS(.4,-20) \arl+(0,39),
\POS(75,-21) \arl+(0,3), \POS(100,-21) \arl+(0,3),  \POS(-1,0) \arl+(3,0), 
\POS(75,10) \arld-(76,0) \arld-(0,10), \POS(100,0) \arld-(25,0) \arld-(0,8),
\POS(100,-8) \arlp+(15,0)
\endxy
    \label{4.3.1} \end{equation}

	We use again the strictly increasing sequence $ x_n = (n+1)/(n+2)$ in $X$, which begins at 
$ x_0 = 1/2 $ and converges to $1$ in $ \bbI$. We choose a decreasing sequence $ \ep = (\ep_n) $ 
such that 
$$
0  <  \ep_n  <  1/2,   \qq   w_n  =  (x_n, \ep_n)  \in  W_n,
$$
and we let $ W'_n = [0, x_n[ \times [0, \ep_n[$. Thus $W'_\ep = \Cup \, W'_n$ is an open nbd of $ X \ti \{0\} $ 
in $ X \ti Y$, and every $ W_n $ has a point $ w_n \notin W'_\ep$.

	Finally $ W'_\ep \cup ([0, 1/2[ \times Y) $ is an open nbd of $ X \jo Y $ in $ X \ti Y $ and does not 
contain any $ W_n$.	
\end{proof}
%

\subsection{A non-countable local basis}\label{4.4}
Letting $X = Y = [0, 1[$, for simplicity, we prove now that the following family of open sets, constructed 
with the same sequence $x_n = (n+1)/(n+2)$ and indexed by the continuum $E$ of all (weakly) decreasing 
sequences $ \ep = (\ep_n) $ of the interval $]0, 1/2[$
    \begin{equation} \begin{array}{l}
W_\ep  =  W'_\ep \cup W''_\ep
\\[5pt]
\;\;\;\;\;\,   = \Cup_n \, (([0, x_n[ \times [0, \ep_n[ ) \cup ([0, \ep_n[ \times [0, x_n[ ))  \q   (\ep \in E),
    \label{4.4.1} \end{array} \end{equation}
%
%
    \begin{equation*} 
\xy <.4mm, 0mm>:
%
(-15,30) *{}; (0,-32) *{};
(30,-25) *{\sst{\ep_0}}; (75,-25) *{\sst{x_0}}; (78,15) *{\sst{w_0}}; (75,10) *{\bu}; 
(100,-25) *{\sst{x_1}}; (103,5) *{\sst{w_1}}; (100,0) *{\bu}; 
(-7,10) *{\sst{\ep_0}}; (-7,0) *{\sst{\ep_1}}; 
(170,-15) *{X \ti X}; (18,-2) *{W_\ep};
\POS(0,-20) \arl+(149,0) \arl+(0,44), \POS(150,-20) \arld+(0,45), 
\POS(0,-19.6) \arl+(149,0),  \POS(.4,-20) \arl+(0,44),
\POS(30,-21) \arl+(0,3), \POS(75,-21) \arl+(0,3), \POS(100,-21) \arl+(0,3),  \POS(-1,0) \arl+(3,0), 
\POS(30,10) \arld+(0,15), \POS(75,10) \arld-(45,0) \arld-(0,10), \POS(100,0) \arld-(25,0) \arld-(0,8),
\POS(100,-8) \arlp+(15,0)
\endxy
    \label{4.4.1bis} \end{equation*}
forms a basis of open nbds of $ X \jo X $ in $ X \ti X$, so that their projections in $ X \me X $ form a 
local basis of the basepoint. We remark that all $ W_\ep$ are symmetric with respect to the diagonal 
of $X^2$.

	In fact, for the compact interval $ X_n = [0, x_n] $ pointed at 0, the open nbds of 
$ X_n \jo X_n $ in $ X_n \ti X_n $ have a basis of L-shaped spaces
    \begin{equation} \begin{array}{l}
L_\ep  =  \{(x, y) \in X_n \ti X_n \sep \; |x| \me |y| < \ep\}
\\[5pt]
\;\;\;\;\;   =  ([0, x_n] \times [0, \ep[ ) \cup ([0, \ep[ \times [0, x_n])   \q   (0 < \ep \le x_n).
    \label{4.4.2} \end{array} \end{equation}

	Now, given any nbd $ W $ of $ X \jo X $ in $ X \ti X$, there is some $\ep_n$ such that
    \begin{equation}
L_{\ep_n} \sub  W \cap X_n \ti X_n   \qq   ( \ep_n \in \; ]0, 1/2[),
    \label{4.4.3} \end{equation}
and we can make the sequence $ \ep = (\ep_n) $ decreasing (replacing $ \ep_n $ by 
$\min_{k\le_n} \ep_k$). Now the union of the interior parts of all $ L_{\ep_n}$ gives a nbd 
$ W_\ep $ of our system, contained in $ W$
    \begin{equation} \begin{array}{c}
W_\ep  \, = \,  \Cup_n \,  \int(L_{\ep_n})  \sub  W,			
\\[5pt]
\int(L_{\ep_n})  =  ([0, x_n[ \times [0, \ep_n[ ) \cup ([0, \ep_n[ \times [0, x_n[ ).
    \label{4.4.4} \end{array} \end{equation}

\LL \begin{small}

\vskip -2pt

\ndt {\em Note.} The local basis $ (L_\ep) $ in \eqref{4.4.2} can be made countable, requiring $ \ep $ 
to be rational, but this is of no relevance here: the set of all decreasing sequences $(\ep_n)$ of 
rational numbers in  ]0, 1/2[  is still a continuum.

\end{small} \LB

\subsection{A larger local basis}\label{4.5}
One can construct a larger basis of open nbds of $ X \jo X $ in $ X \ti X$, which is easily described in 
the new variables $u, v$ (by a $45^{\circ}$ rotation and homothety)
%
%
    \begin{equation} 
\xy <.4mm, 0mm>:
%
(-45,45) *{}; (0,-45) *{};
(-30,-36) *{\sst{-1}}; (30,-36) *{\sst{1}}; 
(-38,1) *{\sst{y}}; (38,1) *{\sst{x}}; (42,-36) *{\sst{u}}; (5,34) *{\sst{v}}; 
(100,15) *{u  =  x - y}; (100,0) *{v  =  x + y}; (34,-18) *{X \ti X}; 
(0,-15); (-30,0) **\crv{~*=<5pt>{.} (-10,-15)&(-20,5)},    
(0,-15); (30,0)  **\crv{~*=<5pt>{.} (10,-15)&(20,5)},
\POS(0,-30) \arl+(-29,29) \arl+(29,29), \POS(0,30) \arld+(-30,-30) \arld+(30,-30), 
\POS(0,-29.6) \arl+(-29,29) \arl+(29,29),
\POS(0,-30) \ard+(-40,40) \ard+(40,40), 
\POS(-42,-30) \ard+(92,0), \POS(0,-38) \ard+(0,78), 
\POS(-30,-31) \arl+(0,2), \POS(30,-31) \arl+(0,2),  
\endxy
    \label{4.5.1} \end{equation}
$$
X \ti X  =  \{(u, v) \sep  |u| < 1,  |u| \le v < 2 - |u|\}.
$$

	Every continuous mapping $ v = f(u) $ between $ v = |u| $ and $ v = 2 - |u| $
    \begin{equation}
f\c ]\! - 1, 1[ \to \bbR,    \qq    |u|  \le  f(u)  \le  2 - |u|,
    \label{4.5.2} \end{equation}
produces an open set $ W(f) $ of $ X \ti X $ containing $ X \jo X$
    \begin{equation} \begin{array}{l}
W(f)	=  \{(u, v) \in \bbR^2 \sep  |u| < 1, |u| \le v < f(u)\}  
\\[2pt]
\q\;   =  \{(x, y) \in X \ti X \sep  x + y < f(x - y)\}.
    \label{4.5.3} \end{array} \end{equation}

	This family contains the previous family $ (W_\ep) $ of \eqref{4.4.4}, taking into account that the 
boundary of $ W_\ep $ in $ X \ti X $ is the graph of a map $ v = f_\ep(u) $ (in the new variables), 
piecewise linear on each compact subinterval of $\, ]\! - 1, 1[$. (We can require that $ f $ be an even 
function, like all functions $ f_\ep$.)

\vskip 3pt

	We conclude this part extending Theorem \ref{4.2}.

\subsection{Theorem {\rm (Non-exponentiable smash products)}}\label{4.6}
{\em
Let $X$ and $\, Y$ be pointed spaces and suppose that:

\LL

\vskip -2pt

\ndt (i) $ X $ has an infinite closed discrete subspace $D$,

\Ndt (ii) $Y$ has a discrete subspace $E$ and $\overline{E} = E \cup \{0_Y\}$ is a disjoint union; 
moreover every point of $E$ is closed in $\overline{E}$.

\LB

\skp

	Then $X \me Y $ is not locally compact: the basepoint $[0]$ has no local basis of compact nbds.

\vskip 3pt

	If, moreover, $ X $ and $ Y $ are Hausdorff spaces, $ X \me Y $ is not exponentiable.
}
\begin{proof}
	As in \ref{4.2}, $p\c X \ti Y \to X \me Y$ is the canonical projection. An arbitrary nbd of $ [0] $ 
in $X \me Y$ is the image $ p(W) $ of a nbd $ W $ of $ X \jo Y $ in $ X \ti Y$,  and we prove that 
$p(W)$ cannot be compact, constructing a closed subset $p(C)$ which is not compact.

	We choose an injective sequence $ (x_n) $ in $ D \setm \{0_X\}$; we can assume that $ D $ 
is the support of this sequence, which is still an infinite closed discrete subspace of $ X$.

	For every $ n \ge 0$, $W \cap (\{x_n\} \ti Y) $ is a nbd of $ (x_n, 0) $ in $ \{x_n\} \ti Y$;  we choose 
a point $ y_n \in E $, and form the subset $C$
    \begin{equation}
w_n  =  (x_n, y_n)  \in  W \setm (X \jo Y),   \q\;\;   C = \{w_n \sep n \ge 0\} \sub D \ti E.
    \label{4.6.1} \end{equation}

	The space $C$ is discrete: every point $(x_n, y) \in D \ti Y$ has a nbd $U_n \ti Y$ in $X \ti Y$ that 
only meets $C$ at $(x_n, y_n)$.
It is also closed in $D \ti \overline{E}$ (and $ X \ti Y$): every point $(x_n, y)$ with $y \neq y_n$ has a nbd 
$U_n \times (\overline{E} \setm \{y_n\})$ in $X \ti \overline{E}$ that does not meet $C$ (because all $y_n$ 
are closed in $\overline{E}$).

	Now the subset $ p(C) $ is closed in $ X \me Y $ (because $ C $ does not meet $ X \jo Y$) and 
in $ p(W)$. The compactness of the latter would imply the compactness of $ p(C)$, which is 
discrete and infinite (because the projection $p$ is injective on $C$).

	The last statement is an obvious consequence.
\end{proof}
%

\subsection{Other examples}\label{4.7}
More explicitly, we want to describe a class of {\em non-associative} triples of pointed spaces of the form 
$(X, Y, \bbQ)$ (or $(\bbQ, Y, X)$, by symmetry).

	We assume that $X$ and $Y$ satisfy the following conditions, similar to 
those of the previous theorem:

\LL
\vskip -2pt

\ndt (i) $X$ has an infinite closed discrete subspace $D$,

\vskip 2pt

\ndt (ii) $Y$ has an infinite discrete subspace $E$ and $\overline{E} = E \cup \{0\}$ is a first-countable 
subspace with $0 \notin E$.

\LB
\vskip 5pt

	For instance these conditions are satisfied when:

\LL
\vskip -2pt

\ndt - $ X $ is $ \bbN$, $\bbZ$, $\bbQ$, $\bbR $ or any {\em non-compact} real interval pointed at any 
point,

\ndt - $ Y $ is $ \bbQ$, $\bbI$, $\bbR$, $\bbS^n$ $ (n > 0) $ or any non-degenerate real interval pointed 
at any point.

\LB
\vskip 5pt

	(Each of these lists can be closed under finite non-empty products.)
	
\vskip 2pt

	On the other hand, we have seen in Theorem \ref{2.2}(b) that any triple $(X', Y', Z) $ of pointed 
spaces where $ X', Y' $ are compact and exponentiable (e.g.\ compact Hausdorff) is regularly associative.

\subsection{Theorem}\label{4.8}
{\em
Let $ X $ and $ Y $ be pointed spaces verifying the previous conditions (i), (ii), respectively.

\Ndt (a) The functor $ - \times \bbQ $ does not preserve the quotient $\, p\c X \ti Y \to X \me Y$, and the smash 
product $ (X \me Y) \me \bbQ $ is not regular.

\Ndt (b) If moreover $ X $ is exponentiable, the triple $ (X, Y, \bbQ) $ is not associative. By symmetry, 
$ (\bbQ, Y, X) $ is neither, including the well-known case $(\bbQ, \bbQ, \bbN)$.
}
\begin{proof}
(a) We choose an injective sequence $(x_n)$ in $D$ which does not contain 0, and an injective 
sequence $(y_n)$ in $E$ that converges to 0 in $\overline{E}$ (and $y$), noting that all $x_n$ are open 
in $D$ and all $y_n$ are open in $E$. We also choose a sequence $(r_n)$ of rational numbers converging 
to $\sqrt{2}$.

	The argument will be based on the following maps (see \eqref{1.8.2})
    \begin{equation} \begin{array}{c}
p \ti 1\c X \ti Y \ti \bbQ \to (X \me Y) \ti \bbQ, 
\\[5pt]
\rho = q(p \ti 1)\c X \ti Y \ti \bbQ \to (X \me Y) \me \bbQ,
    \label{4.8.1} \end{array} \end{equation}
and the sets
    \begin{equation} \begin{array}{c}
C  =  \{(x_m, y_n, r_n/m) \sep  m, n > 0\}  \sub  D \ti E \ti \bbQ,
\\[5pt]
C'  =  (p \ti 1)(C)  \sub  (X \me Y) \ti \bbQ.
    \label{4.8.2} \end{array} \end{equation}

	We note that $C$ does not meet $X \jo Y \jo \bbQ$ and is saturated for the previous maps.

\Ndt (i) It is easy to see that $C$ is closed in $D \ti \overline{E} \ti \bbQ$, a first-countable subspace of 
$X \ti Y \ti \bbR$, and therefore is closed in $X \ti Y \ti \bbQ$.

	In fact, for a sequence $((x_{m_k}, y_{n_k}, r_{n_k}/{m_k}))_k$ of $C$ that converges in $X \ti Y \ti \bbR$, 
the sequence $(m_k)$ is eventually constant at some $m$, while the sequence $(n_k)$ is eventually 
constant at some $n$ or tends to $\infty$.

	Therefore the given sequence converges to $(x_m, y_n, r_n/m)$, which belongs to $C$, or to 
$ (x_m, 0, \sqrt{2}/m)$, which does not belong to $D \ti \overline{E} \ti \bbQ$.

\Ndt  (ii) We prove that $(p \times 1)(0, 0, 0)$ belongs to the closure of $C'$ in $(X \me Y) \ti \bbQ$,  
so that $C'$ is not closed there.

\skp

	Let $N'$ be a nbd of $(0, 0)$ in $(X \me Y) \ti \bbQ$ and $N = (p \ti 1)^{-1}(N')$  the corresponding 
nbd of $0 = (0, 0, 0)$ in $X \ti Y \ti \bbQ$, saturated for $p \ti 1$. There is a basic nbd $U \ti V \ti W$  of  0  
contained in $N$, and we fix $m > 0$ so that $r_n/m \in W$ for all $n > 0$ (as the sequence $(r_n)$ 
spans a bounded subset of $ \bbQ$). Thus $(0, 0, r_n/m) \in N$ for all $n$. By saturation 
$(x_m, 0, r_n/m) \in N$; we can find a basic nbd $U_n \ti V_n \ti W_n$ of each such point in $N$, 
and we can assume that $W_n \sups W_{n'}$ for $n' \ge n$.

\skp
\begin{small}

	For the last point we use a nested local basis of $\sqrt{2}/m$ in $\bbR$ and its traces $W_n$ on 
$\bbQ$; replacing $(r_n)$ by a subsequence we get $r_n/m \in W_n$, for $n > 0$.

\end{small} \LB

\skp

	We choose a pair $k \ge n$ so that $y_k \in V_n$. Finally $(x_m, y_k, r_k/m)$ belongs to 
$U_n \ti V_n \ti W_k \sub U_n \ti V_n \ti W_n \sub N$,  and its projection by $p \ti 1$ belongs to $N'$, 
which meets $C'$.

\Ndt (iii) The argument is done: $C'$ is not closed in $(X \me Y) \ti \bbQ$ and $p \ti 1$ is not a quotient 
map. But $C$ is also saturated for $\rho$, which is not a quotient map either (by an equivalent condition 
in \ref{1.8}(b)), showing that the smash product $(X \me Y) \me \bbQ$ is not regular.

\Ndt (b) A consequence, because $ X \me (Y \me \bbQ) $ is regular while $ (X \me Y) \me \bbQ $ is not.
\end{proof}

\vskip 10pt


\skp

\end{document}